\RequirePackage[table]{xcolor}
\documentclass{siamart0516}
\usepackage{geometry}
\usepackage[table]{xcolor}
\usepackage{graphicx}
\usepackage{amssymb}
\usepackage{amsmath}
\usepackage{bm}
\usepackage{caption}
\usepackage{subcaption}
\usepackage{booktabs}
\usepackage{tikz}
\usepackage{tikzscale}
\usepackage{xspace}
\usepackage{pgfplots}
\usepackage{multirow}
\usepackage{framed}

\usepackage{hyperref}

\usetikzlibrary{shapes.geometric}
\usetikzlibrary{positioning}
\usepgfplotslibrary{external} 

\newcommand\nth{th\xspace} 

\newcommand{\R}{\mathbb{R}}

\def\arraystretch{1.75}

\title{On the convergence of iterative solvers for polygonal discontinuous 
Galerkin discretizations}
\author{Will Pazner, Per-Olof Persson}
\date{}

\begin{document}

\maketitle
\begin{abstract}
We study the convergence of iterative linear solvers for
discontinuous Galerkin discretizations of systems of hyperbolic conservation
laws with polygonal mesh elements compared with that of traditional triangular
elements. We solve the semi-discrete system of equations by means of an
implicit time discretization method, using iterative solvers such as the block
Jacobi method and GMRES. We perform a von Neumann analysis to analytically study
the convergence of the block Jacobi method for the two-dimensional advection
equation on four classes of regular meshes: hexagonal, square,
equilateral-triangular, and right-triangular. We find that hexagonal and square
meshes give rise to smaller eigenvalues, and thus result in faster convergence
of Jacobi's method. We perform numerical experiments with variable velocity
fields, irregular, unstructured meshes, and the Euler equations of gas dynamics
to confirm and extend these results. We additionally study the effect of
polygonal meshes on the performance of block ILU(0) and Jacobi preconditioners
for the GMRES method.
\end{abstract}

\section{Introduction}

In recent years, the Discontinuous Galerkin (DG) method has become a popular
choice for the discretization of a wide range of partial differential equations
\cite{Reed_Hill,cockburn01rkdg,hesthaven08dgbook}. This is partly because of its
many attractive properties, such as the arbitrarily high degrees of
approximation, the rigorous theoretical foundation, and the ability to use fully
unstructured meshes. Also, due to its natural stabilization mechanism based on
approximate Riemann solvers, it has in particular become widely used in fluid
dynamics applications where the high-order accuracy is believed to produce
improved accuracy for many problems \cite{Wang2013}.

Most work on DG methods has been based on meshes of either simplex elements
(triangles and tetrahedra), block elements (quadrilaterals and hexahedra), or
combinations of these such as prism elements. This is likely because of the
availability of excellent automatic unstructured mesh generators, at least for
the simplex case
\cite{peraire87advancingfront,ruppert95delaunay,shewchuk02delaunay}, and also
because of the advantages with the outer-product structure of block elements.
However, it is well known that since no continuity is enforced between the
elements, it is straightforward to apply the DG methods to meshes with elements
of any shapes (even non-conforming ones). For example, vertex-centered DG
methods based on the polygonal dual meshes were studied in
\cite{berggren06vertex,luo08taylor}. This is a major advantage over standard
continuous FEM methods, which need significant developments for the extension to
arbitrary polygonal and polyhedral elements \cite{manzini2014polygonal}.

In the finite volume CFD community, there has recently been considerable
interest in meshes of arbitrary polygonal and polyhedral elements. In fact, the
popular vertex-centered finite volume method applied to a tetrahedral mesh can
be seen as a cell-centered method on the dual polyhedral mesh. Because of this,
a number of methods have been proposed for generation of polyhedral meshes,
which in many cases have advantages over traditional simplex meshes
\cite{oaks2000polyhedral,garimella2014polyhedral}. Although it is still unclear
exactly what benefits these elements provide, they have been reported to be both
more accurate per degree of freedom and to have better convergence properties in
the numerical solvers than for a corresponding tetrahedral mesh
\cite{peric2004polyhedral,balafas2014polyhedral}. There have also been studies
showing that vertex-centered schemes are preferred over cell-centered
\cite{diskin2010viscous,diskin2011inviscid}, again indicating the benefits of
polyhedral elements.

Inspired by the promising results for polyhedral finite volume method, and the
fact that DG is a natural higher-order extension of these schemes, in this work
we study some of the properties of DG discretizations on polygonal meshes. To
limit the scope, we only investigate the convergence properties of iterative
solvers for the discrete systems, assuming an equal number of degrees of freedom
per unit area for all element shapes. Future work will also investigate the
accuracy of the solutions on the different meshes. We first consider the
iterative block-Jacobi method applied to a pure convection problem, which in the
constant coefficient case can be solved analytically using von Neumann analysis.
Next we apply the solver to Euler's equations of gas dynamics for relevant model
flow problems, to obtain numerical results for the convergence of the various
element shapes. We consider regular meshes of hexagons, squares, and two
different configurations of triangles, as well as the dual of fully unstructured
triangular Delaunay refinement meshes. We also perform numerical experiments
with the GMRES Krylov subspace solver and a block-ILU preconditioner. Although
the results are not entirely conclusive, most of the results indicate a clear
benefit with the hexagonal and quadrilateral elements over the triangular ones.

The paper is organized as follows. In Section 2, we describe the spatial and the
temporal discretizations, and introduce the iterative solvers. In Section 3 we
perform the von Neumann analysis of the constant coefficient advection problem,
in 1D and for several mesh configurations in 2D. In Section 4 we show numerical
results for more general advection fields, for more general meshes, as well as
for the Euler equations and the GMRES solver. We conclude with a summary of our
findings as well as directions for future work.

\section{Numerical methods}
\subsection{The discontinuous Galerkin formulation}

We consider a system of $m$ hyperbolic conservation laws given by the equation
\begin{equation}
    \label{eq:cons-law}
    \begin{cases}
    \partial_t \bm{u} + \nabla \cdot \bm{F}(\bm{u}) = 0,
        \qquad(t,\bm{x}) \in [0, T] \times \Omega\\
    \bm{u}(0, \bm{x}) = \bm{u}_0(\bm{x}).
    \end{cases}
\end{equation}
In order to describe the discontinuous Galerkin spatial discretization, we
divide the spatial domain $\Omega \subseteq \R^2$ into a collection of
elements, to form the \emph{triangulation} $\mathcal{T}_h = \{ K_i \}$. Often the
elements $K_i$ are considered to be triangles or quadrilaterals, but in this paper
we allow the elements to be arbitrary polygons in order to study the impact
of different tessellations on the efficiency of the algorithm.

Let $V_h = \left\{ v_h \in L_2(\Omega) : v_h \big|_{K_i} \in P^p(K_i) \right\}$ 
denote the space of piecewise polynomials of degree $p$. We let $\bm{V}_h^m$ 
denote the space of vector-valued functions of length $m$, with each component 
in $V_h$. Note that continuity is not enforced between the elements.
We derive the discontinuous Galerkin method by replacing $\bm{u}$ in equation 
\eqref{eq:cons-law} by an approximate solution $\bm{u}_h \in \bm{V}_h^m$, and then 
multiplying equation  by a test function $\bm{v}_h \in \bm{V}_h^m$. We then 
integrate by parts over each element. Because the approximate solution $\bm{u}_h$ 
is potentially discontinuous at the boundary of an element, the flux function 
$\bm{F}$ is approximated by a \emph{numerical flux function} $\widehat{\bm{F}}$, 
which takes as arguments $\bm{u}^+$, $\bm{u}^-$, and $\bm{n}$, denoting the solution on 
the exterior and interior of the element, and the outward-pointing normal vector, 
respectively. Then, the discontinuous 
Galerkin method reads:\\
\indent
Find $\bm{u}_h \in \bm{V}_h^m$ such that, for all $\bm{v}_h \in \bm{V}_h^m$,
\begin{gather} 
    \label{eq:semi-disc-dg}
    \int_{K_i} \partial_t \bm{u}_h \cdot \bm{v}_h ~dx -
        \int_{K_i} \bm{F}(\bm{u}_h) : \nabla \bm{v}_h~dx
        + \oint_{\partial K_i} \widehat{\bm{F}}(\bm{u}^+, \bm{u}^-, \bm{n}) \cdot \bm{v}_h~ds = 0.
\end{gather}

\subsection{Advection equation}
As a first example, we consider the two-dimensional scalar advection equation
\begin{equation}
    \label{eq:2d-advection}
    u_t + \nabla \cdot \left( \bm{\beta} u \right) = 0,
\end{equation}
for a given (constant) velocity vector $\bm{\beta} = (\alpha, \beta)$. We solve 
this equation in the domain $[0, 2\pi] \times [0, 2\pi]$, with periodic boundary 
conditions. The exact solution to this equation is given by
\begin{equation}
    u(t, x, y) = u_0(x - \alpha t, y - \beta t),
\end{equation}
where $u_0$ is the given initial state.

In order to define the discontinuous Galerkin method for equation \eqref{eq:2d-advection}, 
we define the \emph{upwind flux} by
\begin{equation}
    \widehat{\bm{F}}(\bm{u}^+, \bm{u}^-, \bm{n}) =
    \begin{cases}
        \bm{u}^- \quad\text{if $\bm{\beta}\cdot\bm{n} \geq 0$}\\
        \bm{u}^+ \quad\text{if $\bm{\beta}\cdot\bm{n} < 0$}\\
    \end{cases}
\end{equation}
We represent the approximate solution function $\bm{u}_h$ as a vector $\bm{U}$
consisting of the coefficients of the expansion of $\bm{u}_h$ in terms of an
orthogonal Legendre polynomial modal basis of the function space
$\bm{V}_h^m$. Discretizing equation \eqref{eq:2d-advection} results in a linear
system of equations, which we can write as
\begin{equation}
    \label{eq:adv-lin-system}
    \mathbf{M}(\partial_t\bm{U}) + \mathbf{L}\bm{U} = 0,
\end{equation}
where the mass matrix $\mathbf{M}$ corresponds to the first term on the left-hand 
side of \eqref{eq:semi-disc-dg}, and $\mathbf{L}$ consists of the second two 
terms on the left-hand side. The mass matrix is block-diagonal, and the matrix 
$\mathbf{L}$ is a block matrix, with blocks along the diagonal, and off-diagonal 
blocks corresponding to the boundary terms from the neighboring elements.

\subsection{Temporal integration and linear solvers}
We consider the solution of \eqref{eq:adv-lin-system} by means of implicit time 
integration schemes, the simplest of which is the standard backward Euler scheme,
\begin{equation} 
    \label{eq:be}
    (\mathbf{M} + k \mathbf{L})\bm{U}^{n+1} =  \mathbf{M}\bm{U}^n.
\end{equation}
Furthermore, each stage of a higher-order scheme, such as a diagonally-implicit 
Runge-Kutta (DIRK) scheme \cite{Alexander1977}, can be written as a 
similar equation. The block sparse system can be solved efficiently by means of
an iterative linear solver. In this paper, we consider two solvers: the 
simple block Jacobi method, and the preconditioned GMRES method.

\subsubsection{Block Jacobi method}
A popular and simple iterative solver is the block Jacobi method, defined as
follows. Each iteration of the method for solving the linear system
$\mathbf{A}\bm{x} = \bm{b}$ is given by
\begin{equation}
    \label{eq:jacobi}
    \bm{x}^{(n+1)} = \mathbf{D}^{-1}\bm{b} + \mathbf{R}_J\bm{x}^{(n)},
\end{equation}
where $\mathbf{D}$ is the block-diagonal part of $\mathbf{A}$, and $\mathbf{R}_J
= \mathbf{I} - \mathbf{D}^{-1}\mathbf{A}$. This simple method has the advantage
that it is possible to analyze the convergence properties of the method simply
by examining the eigenvalues of the matrix $\mathbf{R}_J$. An upper bound of 1 
for the absolute value of the eigenvalues of the matrix $\mathbf{R}_J$ is a 
necessary and sufficient condition in order for Jacobi's method to converge 
(for any choice of initial vector $\bm{x}^{(0)}$). The spectral radius of
$\mathbf{R}_J$ determines the speed of convergence.

\subsubsection{Preconditioned GMRES method}
Another popular and oftentimes more efficient \cite{Bassi2000} method for solving
large, sparse linear systems is the GMRES (generalized minimal residual) method
\cite{Persson2009}. As with most Krylov subspace methods, the choice of
preconditioner has great impact on the efficiency of the solver
\cite{Persson-ILU0}. A simple and popular choice of preconditioner is the block
Jacobi preconditioner. Each application of this preconditioner is performed by
multiplying by the inverse of the block-diagonal part of the matrix. Another,
often more effective choice of preconditioner is the block ILU(0) preconditioner
\cite{Diosady2009}. This preconditioner produces an approximate block-wise LU
factorization, whose sparsity pattern is enforced to be the same as that of the
original matrix. This factorization can be performed in-place, and requires no
more storage that the original matrix. Unlike the block Jacobi method, the block
ILU(0) preconditioner can be highly sensitive to the ordering of the mesh
elements \cite{Duff1989,Benzi1999}. Because of this property, it is common to combine the use
of ILU preconditioners with certain orderings of the mesh elements designed to
increase efficiency, such as reverse Cuthill-McKee \cite{Cuthill1969}, minimum
degree \cite{Markowitz1957}, nested dissection \cite{George1973}, or minimum
discarded fill \cite{Persson2009}.

In this paper, we focus our study on the block Jacobi method, which is simpler
and more amenable to analysis. We then perform numerical experiments using both
the block Jacobi method and the preconditioned GMRES method using ILU(0) and
block Jacobi preconditioning.

\section{Jacobi Analysis}

\begin{figure}[h]
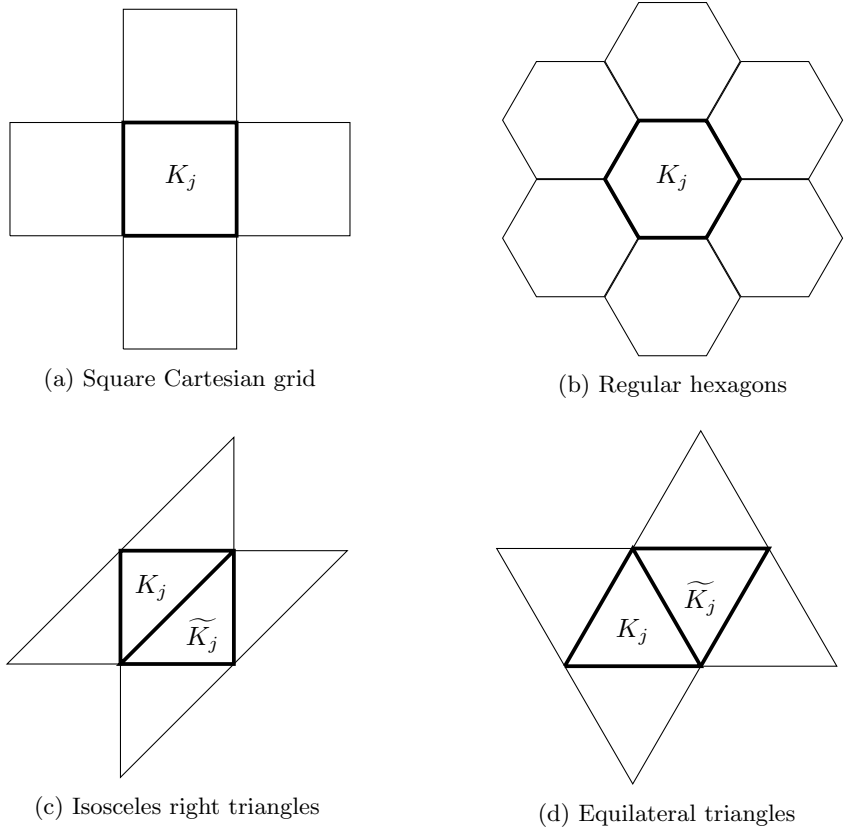

    \centering
    \hspace*{\fill}%
    \begin{subfigure}{0.3\textwidth}
        \centering
        \includegraphics[width=\linewidth]{quad.tikz}
        \caption{Square Cartesian grid}
    \end{subfigure}\hfill%
    \begin{subfigure}{0.3\textwidth}
        \centering
        \includegraphics[width=\linewidth]{hex.tikz}
        \caption{Regular hexagons}
    \end{subfigure}%
    \hspace*{\fill}%
    
    \hspace*{\fill}%
    \begin{subfigure}{0.3\textwidth}
        \centering
        \includegraphics[width=\linewidth]{rtri.tikz}
        \caption{Isosceles right triangles}
    \end{subfigure}\hfill%
    \begin{subfigure}{0.3\textwidth}
        \centering
        \includegraphics[width=\linewidth]{etri.tikz}
        \caption{Equilateral triangles}
    \end{subfigure}
    \hspace*{\fill}%
    
    \caption{Examples of generating patterns $G_j$ shown with bolded lines. 
    Neighboring elements are shown unbolded.}
    \label{fig:gen-pat}
\end{figure}

We compare tessellations of the plane by four sets of \textit{generating
patterns}, each consisting of one or more polygons. We consider tessellations
consisting of squares, regular hexagons, two right triangles, and two
equilateral triangles. The generating patterns considered are shown in Figure
\ref{fig:gen-pat}. Each generating pattern $G_j$ consists of one or two
elements, labeled $K_j$ and $\widetilde{K_j}$. We will refer to these generating
patterns as $S, H, R,$ and $E$ for squares, hexagons, right triangles, and
equilateral triangles, respectively.

We are interested in computing the spectral radius of the Jacobi matrix
$\mathbf{R_J}$ that arises from the discontinuous Galerkin discretization on the
mesh resulting from tessellating the plane by each of the four generating
patterns. For the sake of comparison, we choose the elements from each of the
generating patters to have the same area. Therefore, if the side length of the
equilateral triangle is $h_E = h$, then the two equal sides of the isosceles
right triangle have side length $h_R = \frac{\sqrt[4]3}{\sqrt 2}h_E$, the
hexagon has side length $h_H = \frac{1}{\sqrt 6}h_E$, and the square has side
length $h_S = \frac{\sqrt[4]{3}}{2}h_E$. Then, the global system will have the
same number of degrees of freedom regardless of choice of generating pattern.

\subsection{Von Neumann analysis}
First, we compare the efficiency of each of the four types of generating patterns when
used to solve the advection equation \eqref{eq:2d-advection} with the
discontinuous Galerkin spatial discretization and implicit time integration. We
compute the spectral radius of the matrix $\mathbf{R_J}$ using the classical von
Neumann analysis for each of the generating patterns, in a manner similar to 
\cite{Kubatko2008}.

Let $\bm{U}$ denote the solution vector, and let its $j$\nth component,
$\bm{U}_j$, which is itself a vector, denote the degrees of freedom in $G_j$,
the $j$\nth generating pattern. We remark that in the case of squares and
hexagons, this corresponds exactly to the degrees of freedom in the element
$K_j$, but in the case of the triangular generating patters, this corresponds to
the degrees of freedom from both of the elements $K_j$ and $\widetilde{K_j}$. In
order to determine the eigenvalues of $\mathbf{R_J}$, we consider the planar
wave with wavenumber $(n_x, n_y)$ defined by
\begin{equation}
    \bm{U}_j = e^{i(n_x x_j + n_y y_j)} \widehat{\bm{U}},
\end{equation}
where $(x_j, y_j)$ are fixed coordinates in
$G_j$. Then, we let $\ell$ index the
generating patterns neighboring $G_j$, and we let $\bm{\delta}_\ell =
(\delta_{x\ell}, \delta_{y\ell}) = (x_j - x_\ell, y_j - y_\ell)$
be the offsets satisfying $G_j + \bm{\delta}_\ell = G_\ell$. We can then
write the solution in each of the neighboring generating patterns as
\begin{equation}
    \bm{U}_\ell = e^{i(n_x \delta_{x\ell} + n_y \delta_{y\ell})} \bm{U}_j.
\end{equation}
In this case we write the semi-discrete equations \eqref{eq:adv-lin-system} in
the following compact form
\begin{equation}
    \mathbf{M}_j(\partial_t \bm{U}_j) + 
    \sum_\ell e^{i(n_x \delta_{x\ell} + n_y \delta_{y\ell})} \mathbf{L}_{j\ell} \bm{U}_j = 0,
\end{equation}
where the summation over $\ell$ ranges over all neighboring generating patterns,
$\mathbf{M}_j$ denotes the diagonal block of $\mathbf{M}$ 
corresponding to the $j$th generating pattern, and
$\mathbf{L}_{j\ell}$ denotes the block of $\mathbf{L}$ in the $j$\nth row and
$\ell$\nth column. We can write
\begin{equation}
    \widehat{\mathbf{L}}_j = \sum_\ell e^{i(n_x \delta_{x\ell} + n_y \delta_{y\ell})} \mathbf{L}_{j\ell}
\end{equation}
to further simplify and obtain
\begin{equation}
    \label{eq:compact-form}
    \mathbf{M}_j(\partial_t \widehat{\bm{U}}) + \widehat{\mathbf{L}}_j \widehat{\bm{U}} = 0.
\end{equation}
In order to solve equation \eqref{eq:compact-form} using an implicit method, we
consider the backward Euler-type equation
\begin{equation}
    (\mathbf{M}_j + k \widehat{\mathbf{L}}_j) \widehat{\bm{U}}^{n+1} 
        = \mathbf{M}_j \widehat{\bm{U}}^n.
\end{equation}
The Jacobi iteration matrix $\mathbf{R_J}$ can then be written as
\begin{equation}
    \label{eq:compact-jacobi}
    \widehat{\mathbf{R_J}}_j 
        = \mathbf{I} - \mathbf{D}^{-1}(\mathbf{M}_j + k \widehat{\mathbf{L}}_j),
\end{equation}
where the matrix $\mathbf{D} = \mathbf{M}_j + k \mathbf{L}_{jj}$ consists of the
$j$\nth diagonal block of $\mathbf{M} + k \mathbf{L}$. The eigenvalues of the
matrix $\widehat{\mathbf{R_J}}_j$ control the speed of convergence of Jacobi's
method. In the simple cases of piecewise constant functions ($p=0$), or in the
case of a one-dimensional domain, the eigenvalues can be computed explicitly. In
the more complicated case of $p \geq 1$ in a two-dimensional domain, we compute
the eigenvalues numerically.

\subsection{1D example}\label{sec:1d}
To illustrate the von Neumann analysis, we consider the one-dimensional scalar
advection equation
\begin{equation}
    u_t + u_x = 0
\end{equation}
on the interval $[0, 2\pi]$ with periodic boundary conditions. We divide the
domain into $N$ subintervals $K_j$, each of length $h$. Let $\bm{U}$ denote the
solution vector, and let $\bm{U}_j$ denote the degrees of freedom for the
$j$\nth interval $K_j$. For example, if piecewise constants are used, the method
is identical to the upwind finite volume method, and each $\bm{U}_j$ represents
the average of the solution over the interval. If piecewise polynomials of
degree $p$ are used, each $\bm{U}_j$ is a vector of length $p+1$.

For the purposes of illustration, we choose $p=1$, and let $\bm{U}_j = (u_{j,1},
u_{j,2})$ represent the value of the solution at the left and right endpoints of
the interval $K_j$. Then, the local basis on the interval $K_j$ consists of the
functions 
\begin{equation}
    \label{eq:1d-basis}
    \phi_{j,1}(x) = j - x/h, \hspace{1in} \phi_{j,2}(x) = x/h - j + 1.
\end{equation}
We remark that the upwind flux in this case is always equal to the value of the
function immediately to the left of the boundary point:
\begin{equation}
    \left[ \widehat{\bm{F}}(u^+, u^-, x) v(x) \right]_{(j-1)h}^{jh}
        = u_{j,2}v_{j,2} - u_{j-1,2}v_{j,1}.
\end{equation}

The entries of the $j$\nth block of the mass matrix $\mathbf{M}$ are given by
\begin{equation}
    (\mathbf{M}_j)_{i\ell} = \int_{(j-1)h}^{jh} \phi_{j,i}(x) \phi_{j,\ell}(x)~dx.
\end{equation}
Additionally, we remark that the diagonal blocks of $\mathbf{L}$ consist of the 
volume integrals and right boundary terms given by
\begin{equation}
    (\mathbf{L}_{jj})_{i\ell} = \phi_{j,i}(jh)\phi_{i,\ell}(jh) 
        - \int_{(j-1)h}^{jh} \phi'_{j,i}(x)\phi_{j,\ell}(x)~dx.
\end{equation}
We let $\bm{A}$ denote the backward Euler-type operator defined by
\begin{equation}
    \mathbf{A} = \mathbf{M} + k\mathbf{L},
\end{equation}
and, solving the equation $\mathbf{A}\bm{x} = \bm{b}$ by means of Jacobi iterations, 
we define the Jacobi matrix $\mathbf{R_J}$ by
\begin{equation}
    \mathbf{R_J} = \mathbf{I} - \mathbf{D}^{-1}\mathbf{A},
\end{equation}
where $\mathbf{D}$ is the matrix consisting of the diagonal blocks of $\mathbf{A}$. 
The entries of the diagonal blocks $\mathbf{M}_j$ and $\mathbf{L}_{jj}$ can be computed 
explicitly using \eqref{eq:1d-basis} to obtain
\begin{equation}
    \mathbf{M}_j = \left(\begin{array}{cc}
    \frac{h}{3} & \frac{h}{6} \\
    \frac{h}{6} & \frac{h}{3}
    \end{array}\right), \quad
    \mathbf{L}_{jj} = \left(
    \begin{array}{cc}
     \frac{1}{2} & \frac{1}{2} \\
     -\frac{1}{2} & \frac{1}{2} \\
    \end{array}
    \right), \quad
    \mathbf{D}_j = \left(\begin{array}{cc}
    \frac{h}{3} + \frac{k}{2} & \frac{h}{6} + \frac{k}{2} \\
    \frac{h}{6} - \frac{k}{2} & \frac{h}{3} + \frac{k}{2}
    \end{array}\right).
\end{equation}

In order to perform the von Neumann analysis, we seek solutions of the form 
$\bm{U}_j = e^{inhj} \widehat{\bm{U}}$, which allows us to explicitly compute the 
form of the matrix $\widehat{\mathbf{L}}_j$. Recalling the compact form from \eqref{eq:compact-form}, 
we obtain
\begin{equation}
    \widehat{\mathbf{L}}_j = \left(\begin{array}{cc}
    \frac{1}{2} & \frac{1}{2}-e^{-ihn} \\
    -\frac{1}{2} & \frac{1}{2}
    \end{array}\right).
\end{equation}
Then, the Jacobi matrix $\widehat{\mathbf{R_J}}_j$ is given by
\begin{equation}
    \widehat{\mathbf{R_J}}_j =
    \left(
\begin{array}{cc}
 0 & \frac{2 e^{-i h n} k (2 h+3 k)}{h^2+4 k h+6 k^2} \\
 0 & -\frac{2 e^{-i h n} (h-3 k) k}{h^2+4 k h+6 k^2} \\
\end{array}
\right),
\end{equation}
whose eigenvalues $\lambda_1$ and $\lambda_2$ are given by
\begin{equation}
    \lambda_1 = 0, \qquad \lambda_2 = \frac{2 k (3 k-h) e^{-i h n}}{h^2+4 h k+6 k^2}.
\end{equation}
Therefore, each wavenumber $n$ from $0$ to $2\pi/h$ corresponds to an eigenvalue
of the Jacobi matrix $\mathbf{R_J}$, and the magnitude of these eigenvalues
determine the speed of convergence of Jacobi's method. In this case, the
expression
\begin{equation}
    \lambda_{\rm max} = \frac{2 k \left| h-3 k\right| }{h^2+4 h k+6 k^2}
\end{equation}
determines the speed of convergence of Jacobi's method. This expression can 
easily be seen to be bounded above by 1 for all positive values of $h$ and $k$, 
therefore indicating that Jacobi's method is guaranteed to converge, 
unconditionally, regardless of spatial resolution or timestep. 

\subsection{2D analysis}\label{sec:2d-analysis}
We now turn to the analysis of the four generating patterns shown in Figure 
\ref{fig:gen-pat}. The analysis proceeds along the same lines as in the 
one-dimensional example from Section \ref{sec:1d}. As an example, we present 
the case of piecewise constants, for which it is possible to explicitly compute 
the eigenvalues of the Jacobi matrix $\mathbf{R_J}$. In this case the 
discontinuous Galerkin formulation simplifies to the upwind finite volume method
\begin{equation}
    \int_{K_j} \partial_t u_h~dx  
        + \oint_{\partial K_j} \widehat{\bm{F}}(u^+, u^-, \bm{n})~ds = 0.
\end{equation}
For the sake of concreteness, we assume without loss of generality
that the velocity vector 
$\bm{\beta} = (\alpha, \beta)$ satisfies $\alpha, \beta \geq 0$. In order to 
explicitly write the upwind flux on the meshes 
consisting of hexagons and equilateral triangles, we further assume that 
$\sqrt{3}\alpha - \beta \geq 0$, and on the mesh consisting of right triangles 
we assume that $\alpha - \beta \geq 0$.
In the case of the square and hexagonal meshes, there 
is only one degree of freedom per generating pattern, and we will write 
$u_j$ to represent the average value of the solution over the generating 
pattern $G_j$. We then consider the planar wave with wavenumber 
$(n_x, n_y)$ given by $u_j = e^{i(n_x x_j + n_y y_j)} \widehat{u}$. 
In the case of the square mesh with side length $h_S =
\frac{\sqrt[4]{3}}{2}h_E$, the method can be written as
\begin{equation}
    h_S^2 \left( \partial_t \widehat{u} \right) = 
        - h_S \left( 
            \alpha(1  - e^{-i n_x h_S} )
            + \beta( 1 - e^{-i n_y h_S} ) \right) \widehat{u}.
\end{equation}
In this case, the mass matrix $\mathbf{M}$ is a diagonal matrix with $h_S^2$
along the diagonal, and the diagonal entries of the matrix $\mathbf{L}$ are
given by $h_S(\alpha + \beta)$. Therefore, the eigenvalues of the Jacobi matrix
$\mathbf{R_J^S} = \mathbf{I} - D^{-1}(\mathbf{M} + k\mathbf{L})$ are given by
\begin{equation}
\begin{aligned}
    \label{eq:fv-sq-eigs}
    \lambda(\mathbf{R_J^S}) &= 1 - \frac{1}{h_S^2 + h_Sk(\alpha + \beta)}
        \left(h_S^2 + h_Sk\left( \alpha(1  - e^{-i n_x h_S} )
            + \beta( 1 - e^{-i n_y h_S} ) \right) \right) \\
        &= \frac{k\left( \alpha e^{-in_x h_S} + \beta e^{-in_y h_S} \right)}
                 {h_S + k(\alpha+\beta)}.
\end{aligned}
\end{equation}

In the case of the hexagonal mesh with side length $h_H =
\frac{1}{\sqrt{6}}h_E$, the method is
\begin{equation}
\begin{aligned}
    \frac{3\sqrt{3}}{2}h_H^2 \left(\partial_t \widehat{u} \right) &= 
    - h_H\Bigg(
    \left( \sqrt{3}\alpha + \beta \right)
    +\left( - \tfrac{\sqrt{3}}{2}\alpha + \tfrac{\beta}{2} \right)e^{ih_H\left(-\frac{3}{2}n_x + \frac{\sqrt{3}}{2}n_y\right)}\\
    &\qquad
    +\left( - \tfrac{\sqrt{3}}{2}\alpha - \tfrac{\beta}{2} \right)e^{ih_H\left(-\frac{3}{2}n_x - \frac{\sqrt{3}}{2}n_y\right)}
    -\beta e^{-ih_H\sqrt{3}n_y}
    \Bigg) \widehat{u}.
\end{aligned}
\end{equation}
A similar analysis shows that the eigenvalues of the matrix $\mathbf{R_J^H}$ are 
given by
\begin{equation}
\begin{aligned}
    \label{eq:fv-hex-eigs}
    \lambda(\mathbf{R_J^H})
    &=\tfrac{k e^{-\frac{1}{2} i h_H \left(3 n_x+\sqrt{3} n_y\right)} 
                \left(\sqrt{3} \beta \left(2
                    e^{\frac{1}{2} i h_H \left(3 n_x-\sqrt{3} n_y\right)}
                    - e^{i \sqrt{3} h_H n_y}+1\right)+
                3 \alpha  \left(1+e^{i \sqrt{3} h_H n_y}\right)\right)}
            {9 h_H + 6 \alpha  k+2 \sqrt{3} \beta  k}.
\end{aligned}
\end{equation}

In the case of the two triangular meshes, there are two degrees of freedom per 
generating pattern, corresponding to the elements $K_j$ and $\widetilde{K_j}$ 
in the generating pattern $G_j$. We write $\bm{U}_j = (u_{j,1}, u_{j,2})$, 
where $u_{j,1}$ is the average of the solution over the element $K_j$, and 
$u_{j,2}$ is the average of the solution over $\widetilde{K_j}$. The planar 
wave solution is then given by $\bm{U}_j = e^{i(n_x x_j + n_y y_j)} \widehat{\bm{U}}$, 
for $\widehat{\bm{U}} = (\widehat{u}_1, \widehat{u}_2)$. We consider the case of a 
right-triangular mesh, where the two equal sides of the isosceles right triangles 
have length $h_R = \frac{\sqrt[4]{3}}{\sqrt{2}}h_E$. The method then reads:
\begin{equation}
\partial_t \left(\begin{array}{c} \widehat{u}_1 \\ \widehat{u}_2 \end{array}\right)
= -\frac{2}{h_R}\left(\begin{array}{c}
     \alpha \widehat{u}_1 - e^{-ih_Rn_x} \alpha \widehat{u}_2  \\
     \alpha \widehat{u}_2 + (\beta - \alpha) \widehat{u}_1 
        - e^{-ih_Rn_y} \beta \widehat{u}_1 
\end{array}\right).
\end{equation}
In the case of the mesh consisting of equilateral triangles, each with side
length $h_E$, the method reads:
\begin{equation}
\partial_t \left(\begin{array}{c} \widehat{u}_1 \\ \widehat{u}_2 \end{array}\right)
= \tfrac{-4}{\sqrt{3}h_E}\left(\begin{array}{c}
    \left(\frac{\sqrt 3}{2}\alpha 
        + \frac{1}{2}\beta\right)\widehat{u}_1
    + \left(e^{-ih_En_x}\left(-\frac{\sqrt 3}{2}\alpha + \frac{1}{2}\beta \right)
    - e^{-ih_En_y} \beta\right) \widehat{u}_2 \\
    \left(-\frac{\sqrt 3}{2}\alpha 
    - \frac{1}{2}\beta\right) \widehat{u}_1 + \left(\frac{\sqrt 3}{2}\alpha 
    + \frac{1}{2}\beta \right) \widehat{u}_2
\end{array}\right).
\end{equation}
Computing the eigenvalues of the corresponding Jacobi matrices $\mathbf{R_J^R}$ 
and $\mathbf{R_J^E}$, we obtain
\begin{align}
    \label{eq:fv-rt-eigs}
    \lambda(\mathbf{R_J^R}) &= \pm\frac{2 k e^{-\frac{1}{2} i h_R (n_x+n_y)} 
        \sqrt{\alpha } \sqrt{\beta +(\alpha -\beta ) e^{i h_R n_y}}}{h_R+2 \alpha  k}, \\
    \label{eq:fv-et-eigs}
    \lambda(\mathbf{R_J^E}) &= \pm \frac{2 k \left(3 \alpha +\sqrt{3} \beta \right) \sqrt{2 \beta  e^{i h_E n_x}+\left(\sqrt{3} \alpha -\beta \right) e^{i h_E n_y}}}{\left(3 h_E+6 \alpha  k+2 \sqrt{3} \beta  k\right) \sqrt{\left(\sqrt{3} \alpha +\beta \right) e^{i h_E (n_x+n_y)}}}.
\end{align}

Then, equations \eqref{eq:fv-sq-eigs}, \eqref{eq:fv-hex-eigs},
\eqref{eq:fv-rt-eigs}, and \eqref{eq:fv-et-eigs} completely determine the speed
of convergence for Jacobi's method of each of the four generating patterns
considered. In the case of a higher-order discontinuous Galerkin method with
basis consisting of piecewise polynomials of degree $p>0$, we obtain a Jacobi
matrix given by equation \eqref{eq:compact-jacobi}, where the matrices
$\widehat{\mathbf{R_J}}_j, \mathbf{D}, \mathbf{M}_j,$ and
$\widehat{\mathbf{L}}_j$ are $\frac{(p+1)(p+2)}{2}\times\frac{(p+1)(p+2)}{2}$
blocks. In this case, we do not obtain closed-form expressions for the
eigenvalues, but rather compute them numerically.

\begin{figure}[!b]
   \centering
   \hspace*{\fill}%
   \begin{subfigure}{0.3\textwidth}
       \centering
           \begin{tikzpicture}
           \begin{axis}[
             ymin=0.75,
             ymax=1.35,
             width=1.95in,
             ybar=0pt,
             bar width=7pt,
             xtick=data,
             every axis/.append style={font=\tiny},
             xticklabels={$k_1$, $k_2$, $k_3$},
             enlarge x limits=0.3,
             major x tick style = {opacity=0},
             legend cell align=left,
             legend entries={Hexagons,
                             Squares,
                             Right Triangles,
                             Equilateral Triangles},
             legend columns=4,
             legend style={
                draw=none,
                /tikz/every even column/.append style={column sep=0.5cm}},
             legend to name=leg:barlegend
          ]
          \addplot table[header = false, x index = 0, y index = 1] {deg0.dat};
          \addplot table[header = false, x index = 0, y index = 2] {deg0.dat};
          \addplot table[header = false, x index = 0, y index = 3] {deg0.dat};
          \addplot table[header = false, x index = 0, y index = 4] {deg0.dat};
          \end{axis}
          \end{tikzpicture}
      \caption{$p=0$}
  \end{subfigure}
  \begin{subfigure}{0.3\textwidth}
      \centering
          \begin{tikzpicture}
          \begin{axis}[
             ymin=0.75,
             ymax=1.35,
             width=1.95in,
             ybar=0pt,
             bar width=7pt,
             xtick=data,
             every axis/.append style={font=\tiny},
             xticklabels={$k_1$, $k_2$, $k_3$},
             enlarge x limits=0.3,
             major x tick style = {opacity=0},
          ]
          \addplot table[header = false, x index = 0, y index = 1] {deg1.dat};
          \addplot table[header = false, x index = 0, y index = 2] {deg1.dat};
          \addplot table[header = false, x index = 0, y index = 3] {deg1.dat};
          \addplot table[header = false, x index = 0, y index = 4] {deg1.dat};
          \end{axis}
          \end{tikzpicture}
      \caption{$p=1$}
  \end{subfigure}%
  \hspace*{\fill}%
  
  \hspace*{\fill}%
  \begin{subfigure}{0.3\textwidth}
      \centering
          \begin{tikzpicture}
          \begin{axis}[
             ymin=0.75,
             ymax=1.35,
             width=1.95in,
             ybar=0pt,
             bar width=7pt,
             xtick=data,
             every axis/.append style={font=\tiny},
             xticklabels={$k_1$, $k_2$, $k_3$},
             enlarge x limits=0.3,
             major x tick style = {opacity=0},
          ]
          \addplot table[header = false, x index = 0, y index = 1] {deg2.dat};
          \addplot table[header = false, x index = 0, y index = 2] {deg2.dat};
          \addplot table[header = false, x index = 0, y index = 3] {deg2.dat};
          \addplot table[header = false, x index = 0, y index = 4] {deg2.dat};
          \end{axis}
          \end{tikzpicture}
      \caption{$p=2$}
  \end{subfigure}
  \begin{subfigure}{0.3\textwidth}
      \centering
          \begin{tikzpicture}
          \begin{axis}[
             ymin=0.75,
             ymax=1.35,
             width=1.95in,
             ybar=0pt,
             bar width=7pt,
             xtick=data,
             every axis/.append style={font=\tiny},
             xticklabels={$k_1$, $k_2$, $k_3$},
             enlarge x limits=0.3,
             major x tick style = {opacity=0},
          ]
          \addplot table[header = false, x index = 0, y index = 1] {deg3.dat};
          \addplot table[header = false, x index = 0, y index = 2] {deg3.dat};
          \addplot table[header = false, x index = 0, y index = 3] {deg3.dat};
          \addplot table[header = false, x index = 0, y index = 4] {deg3.dat};
          \end{axis}
          \end{tikzpicture}
      \caption{$p=3$}
  \end{subfigure}
  \hspace*{\fill}%
  
  \tikzexternaldisable\ref{leg:barlegend}\tikzexternalenable
  \caption{Ratios of the logarithm of the largest eigenvalues for 
  each pattern.}
  \label{fig:bars}
\end{figure}
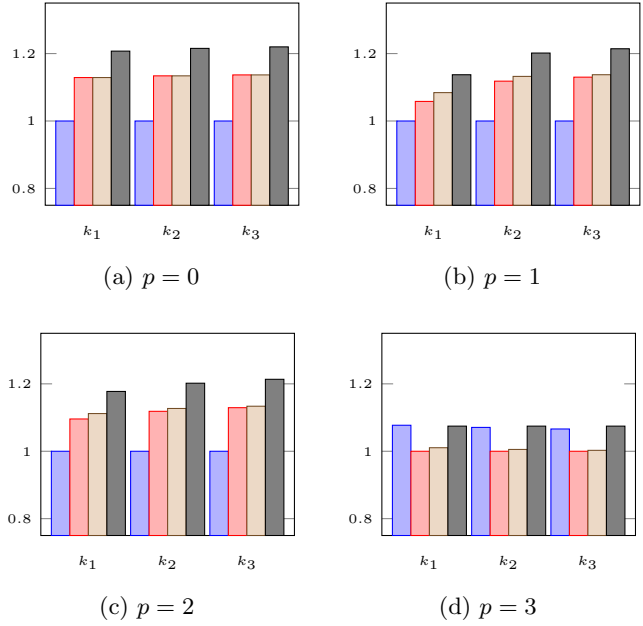

We normalize the velocity magnitude and consider $\bm{\beta} = (\cos(\theta),
\sin(\theta))$. On the square mesh, $\theta$ can range from $0$ to $\pi/2$. On
the hexagonal and equilateral triangle meshes, $\theta$ ranges from $0$ to
$\pi/3$, and on the right-triangular mesh $\theta$ ranges from $0$ to $\pi/4$.
We consider a fixed spatial resolution $h$, and compare the efficiency of the
four patterns for three choices of temporal resolution. We first consider an
``explicit'' time step, satisfying the CFL-type condition
\begin{equation}
    \label{eq:cfl}
    k_{\rm exp} = \frac{h}{|\bm{\beta}|}.
\end{equation}
As one advantage of using an implicit method is that we are not limited by an
explicit timestep restriction of the form \eqref{eq:cfl}, we consider three implicit time
steps given by $k_1 = 3 k_{\rm exp}$, $k_2 = 2 k_1$, and $k_3 = 4 k_1$. We
then maximize over a discrete sample of $\theta \in [0, \pi/4]$ and over all
wavenumbers $(n_x, n_y)$, in order to compute maximum eigenvalue for each of the
generating patterns. As the number of iterations required to converge to a given
tolerance scales like the reciprocal of the logarithm of the spectral radius, we
compare the efficiency of the generating patterns by considering the ratio
\[ \frac{\log\left(\lambda_{\rm max}(\mathbf{R_J^{\mathrm{min}}})\right)}
        {\log\left(\lambda_{\rm max}(\mathbf{R_J^*})\right)}, \]
where $\lambda_{\rm max}(\mathbf{R_J^*})$ is the largest eigenvalue of 
$\mathbf{R_J^*}$, for $* = H, S, R, E$, and 
$\lambda_{\rm max}(\mathbf{R_J^{\mathrm{min}}})$ is the smallest among all 
$\lambda_{\rm max}(\mathbf{R_J^*})$. This ratio corresponds to the ratio of
iterations required to converge to a given tolerance when compared with the most
efficient among the generating patterns. The results obtained for $p=0,1,2,3$,
and $k=k_1,k_2,k_3$ for each generating pattern are shown in Table \ref{tab:eig}
and Figure \ref{fig:bars}.

\newcommand{\highlightcell}{\cellcolor[gray]{0.9}\bf}
\begin{table}[!t]
    \centering
    \tiny
       \begin{tabular}{r | lll | lll}
            & \multicolumn{3}{c|}{\small $p=0$} & \multicolumn{3}{c}{\small $p=1$} \\
            & $k_1$ & $k_2$ & $k_3$ & $k_1$ & $k_2$ & $k_3$ \\
            \hline
            Hexagons    &
\highlightcell1.000000 & \highlightcell1.000000 & \highlightcell1.000000 & 
\highlightcell1.000000 & \highlightcell1.000000 & \highlightcell1.000000\\
            Squares                 & 
1.128939    &   1.133989    &   1.136772    &
1.058098	&	1.118222	&	1.130101\\
            Right triangles         & 
1.128939	&	1.133989	&	1.136772    &
1.084223	&	1.132326	&	1.137313\\
            Equilateral triangles   & 
1.207328	&	1.215467	&	1.219948    &
1.137267	&	1.201638	&	1.214376
        \end{tabular}
       \vspace{8pt}
   \vspace{12pt}
       \begin{tabular}{r | lll | lll}
            & \multicolumn{3}{c|}{\small $p=2$} & \multicolumn{3}{c}{\small $p=3$} \\
            & $k_1$ & $k_2$ & $k_3$ & $k_1$ & $k_2$ & $k_3$ \\
            \hline
            Hexagons                &
\highlightcell1.000000 & \highlightcell1.000000 & \highlightcell1.000000 &
1.077183	&	1.070785	&	1.066101\\
            Squares                 & 
1.095785	&	1.118510	&	1.129314    &
\highlightcell1.000000  &   \highlightcell1.000000  &   \highlightcell1.000000\\
            Right triangles         & 
1.111863	&	1.126951	&	1.133634    &
1.010482	&	1.005391	&	1.002733\\
            Equilateral triangles   & 
1.177503	&	1.201918	&	1.213527    &
1.074570	&	1.074570	&	1.074570
        \end{tabular}
    \vspace{8pt}
   \caption{Ratio of logarithm of eigenvalues
   $\log\left(\lambda_{\rm max}(\mathbf{R_J^{\mathrm{min}}})\right)
   /\log\left(\lambda_{\rm max}(\mathbf{R_J^*})\right)$ ranging 
   over angle $\theta$ and wavenumber $(n_x, n_y)$, for piecewise polynomials 
   of degree 0, 1, 2, and 3, for varying choices of time step $k$. The smallest
   eigenvalue in each column is highlighted.}
    \label{tab:eig}
\end{table}

We remark that for degrees 0, 1, and 2 polynomials, the hexagonal mesh resulted
in the smallest eigenvalues for all choices of timestep considered, and the
square mesh resulted in the second-smallest eigenvalues. For degree 3
polynomials, the square mesh resulted in the smallest eigenvalues for all cases
considered. We notice a significant decrease in the expected performance of the
hexagonal elements in the case of $p = 3$, although we have noticed that the 
effect observed in practice is not as significant as the theoretical results
would suggest.

\pagebreak

\section{Numerical Results}\label{sec:numerical}

\subsection{Advection with variable velocity field} \label{sec:variable-velocity}
\begin{figure}[h]
    \centering
    \includegraphics[width=0.5\textwidth]{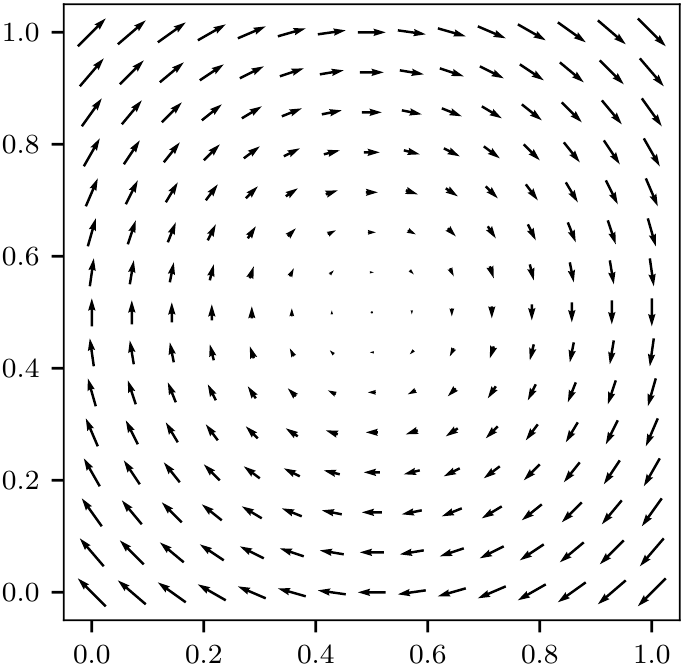}
    \caption{Velocity field $\bm{\beta}(x,y) = (2y - 1, -2x + 1)$}
    \label{fig:vel-field}
\end{figure}

To perform numerical experiments extending the analysis of equation
\eqref{eq:2d-advection} beyond the case of a constant velocity $\bm{\beta}$, we
consider a variable velocity field $\bm{\beta}(x, y)$. In this case, the upwind
numerical flux 
\begin{equation}
    \widehat{\bm{F}}(\bm{u}^+, \bm{u}^-, \bm{n}, x, y) =
    \begin{cases}
        \bm{u}^-(x,y) \quad\text{if $\bm{\beta}(x, y)\cdot\bm{n} \geq 0$}\\
        \bm{u}^+(x,y) \quad\text{if $\bm{\beta}(x, y)\cdot\bm{n} < 0$}\\
    \end{cases}
\end{equation}
is evaluated point-wise. As an example, we define the velocity to be given by the
vector field $\bm{\beta}(x,y) = (2y - 1, -2x + 1)$ on the spatial domain $\Omega
= [0,1]\times[0,1]$. This velocity field is shown in Figure \ref{fig:vel-field}.
We let the initial conditions be given by the Gaussian centered at $(x_0, y_0) =
(0.35, 0.5)$,
\begin{equation}
    u_0(x, y) = \exp(-150((x-x_0)^2 + (y-y_0)^2)).
\end{equation}
The exact solution is periodic with period $\pi$, and is given by the rotation
about the center of the domain,
\begin{equation}
    u(x, y, t) = \exp(-150((x - 0.5 + 0.15\cos(2t))^2 + (y - 0.5 - 0.3\cos(t)\sin(t))^2)).
\end{equation}

\subsubsection{Convergence of the block Jacobi method}
We consider meshes of the domain created by repeating each of the four
generating patterns considered in the previous section. As before, for fixed
spatial resolution $h$, we choose $h_H, h_S, h_R,$ and $h_E$ such that the
number of degrees of freedom is the same for each mesh. We then solve the
advection equation using the backward Euler time discretization, where the block
Jacobi iterative method is used to solve the resulting linear system.
The zero vector is used as the starting vector for the block
Jacobi solver.{ We choose $h = 0.05$, and since $\max_{(x,y)} |\bm{\beta}(x,y)
| = \sqrt{2}$, we consider time steps of $k_1 = h/\sqrt{2}$, $k_2 = 2 k_1$, $k_3
= 4 k_1$. The number of iterations required for the block Jacobi method to
converge to a tolerance of $10^{-14}$ are given in Table
\ref{tab:gaussian-iters}.

\begin{table}[h]
    \centering
    \tiny
    \setlength{\tabcolsep}{8pt}
    \begin{tabular}{r | lll | lll | lll | lll}
    & \multicolumn{3}{c|}{$p=0$} & \multicolumn{3}{c|}{$p=1$} & \multicolumn{3}{c|}{$p=2$}
    & \multicolumn{3}{c}{$p=3$}\\
    & $k_1$ & $k_2$ & $k_3$ & $k_1$ & $k_2$ & $k_3$ & $k_1$ & $k_2$ & $k_3$ & $k_1$ & $k_2$ & $k_3$\\
    \hline
    Hexagons                & \highlightcell33 & \highlightcell57 &
        \highlightcell104 & \highlightcell21 & \highlightcell41 &
        \highlightcell77  & 24 & \highlightcell41 & \highlightcell77 & 
        \highlightcell21 & \highlightcell39 & \highlightcell75 \\
    Squares                 & 35 & 61 & 109 & \highlightcell21 & 42 & 83  
        & \highlightcell22 & 42 & 83  & 22 & 42 & 81 \\
    Right triangles         & 39 & 68 & 128 & 26 & 51 & 100 & 25 & 51 & 100 & 25 & 51 & 100 \\
    Equilateral triangles   & 37 & 67 & 123 & 25 & 47 & 92  & 25 & 47 & 92 & 24 & 47 & 91 \\
    \end{tabular}
    \vspace{12pt}
    \caption{Iterations required for the block Jacobi iterative method to 
        converge in the case of a non-constant velocity field. The smallest
        number of iterations in each column is highlighted.}
    \label{tab:gaussian-iters}
\end{table}

The results are similar to those from the analysis performed in Section
\ref{sec:2d-analysis}. We note that the hexagonal and square meshes resulted in
the lowest number of Jacobi iterations for all of the test cases considered. In
contrast to the results of Section \ref{sec:2d-analysis}, we do not observe a 
decrease in the performance of the hexagonal elements for the case of $p=3$, 
and instead the performance is similar among all choices of $p$ considered.

\subsubsection{Randomly perturbed mesh}
We now consider the effect of polygonal elements on irregular meshes. To this
end, we consider a set of \emph{generating points} distributed evenly on a
Cartesian grid with mesh size $h$. Then, each point is perturbed by a random
perturbation sampled uniformly from the interval $[-\delta, \delta]$. We obtain
two randomized meshes by constructing the Delaunay triangulation and Voronoi
diagram resulting from this set of generating points. The Delaunay mesh consists
entirely of triangular elements, whereas the Voronoi diagram is constructed out
of arbitrary polygonal elements. Examples of the two meshes considered are shown
in Figure \ref{fig:random-mesh}. In contrast to the regular meshes considered in
the previous examples, these two meshes do not consist of the same number of
elements. The Voronoi diagram consists of about half the number of elements as
the Delaunay triangulation. In the test case considered, the randomized
polygonal mesh consists of 410 polygonal elements, whereas the randomized
triangular mesh consists of 759 triangular elements.

The governing equations and set-up is the same as in the previous section. We
record the number of block Jacobi iterations required to converge to  a
tolerance of $10^{-14}$
in Table \ref{tab:random-iters}. Because there is a difference in the number of
mesh elements, the resulting linear system will have a different total number of
degrees of freedom. This difference will then have an additional effect on the
speed of convergence of the block Jacobi method. We note that for polynomials of
degree $p=0,1,2,3$ and for all choices of time step $k$ considered, solving the
system resulting from the Voronoi diagram requires fewer block Jacobi iterations
than does solving the system resulting from the corresponding Delaunay
triangulation.

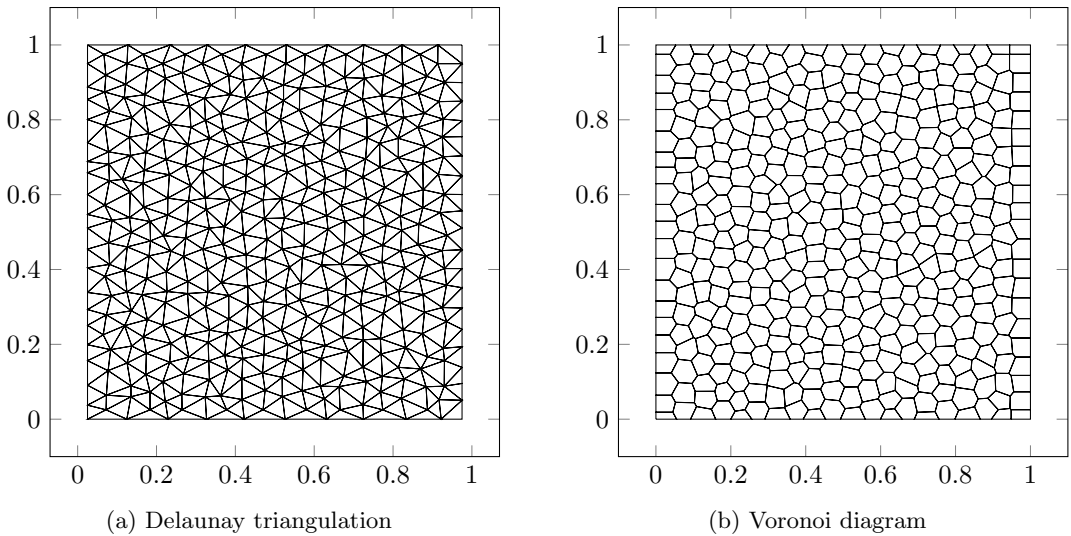
\begin{figure}[h]
    \centering
    \hspace*{\fill}%
    \begin{subfigure}{0.5\textwidth}
        \centering
        \begin{tikzpicture}
        \begin{axis}[
            width=\textwidth,
            height=\textwidth,
            no markers,
        ]
        \addplot[black] table {randomtri.dat};
        \end{axis}
        \end{tikzpicture}
        \caption{Delaunay triangulation}
    \end{subfigure}\hfill%
    \begin{subfigure}{0.5\textwidth}
        \centering
        \begin{tikzpicture}
        \begin{axis}[
            width=\textwidth,
            height=\textwidth,
            no markers,
        ]
        \addplot[black] table {randompoly.dat};
        \end{axis}
        \end{tikzpicture}
        \caption{Voronoi diagram}
    \end{subfigure}%
    \hspace*{\fill}%
    \caption{Randomized polygonal and triangular meshes corresponding to the same 
    set of generating points.}
    \label{fig:random-mesh}
\end{figure}

\begin{table}[h]
    \centering
    \tiny
    \setlength{\tabcolsep}{8pt}
    \begin{tabular}{r | lll | lll | lll | lll}
    & \multicolumn{3}{c|}{$p=0$} & \multicolumn{3}{c|}{$p=1$} & \multicolumn{3}{c|}{$p=2$}
    & \multicolumn{3}{c}{$p=3$}\\
    & $k_1$ & $k_2$ & $k_3$ & $k_1$ & $k_2$ & $k_3$ & $k_1$ & $k_2$ & $k_3$ & $k_1$ & $k_2$ & $k_3$\\
    \hline
    Voronoi diagram &
        \highlightcell27 & \highlightcell32 & \highlightcell38 &
        \highlightcell24 & \highlightcell33 & \highlightcell38  & 
        \highlightcell24 & \highlightcell32 & \highlightcell36 & 
        \highlightcell22 & \highlightcell31 & \highlightcell36 \\
    Delaunay triangulation & 
        38 & 48 & 52 & 
        33 & 45 & 48  &
        33 & 46 & 50 & 
        33 & 44 & 48 \\
    \end{tabular}
    \vspace{12pt}
    \caption{Iterations required for the block Jacobi iterative method to 
        converge in the case of irregular, randomly perturbed meshes. The
        smallest number of iterations in each column is highlighted.}
    \label{tab:random-iters}
\end{table}

\subsubsection{Convergence of the GMRES method} \label{sec:adv-gmres}
The above analysis focused on the block Jacobi method largely because of the
simplicity of the method. In practice, more sophisticated iterative methods are
often used \cite{Persson2009}. In this section, we consider the solution of the
linear system \eqref{eq:be} by means of the GMRES method, using both the block
Jacobi and the block ILU(0) preconditioners. Since the computational work
increases per iteration in GMRES, we choose a \textit{restart parameter} of 20
iterations \cite{saad03iterative}.
We repeat the above test case of the advection equation with variable velocity 
field and record the number of GMRES iterations required to converge to 
a tolerance of $10^{-14}$ using the block Jacobi preconditioner 
in Table \ref{tab:gmres-jac}.

\begin{table}[h]
    \centering
    \tiny
    \setlength{\tabcolsep}{8pt}
    \begin{tabular}{r | lll | lll | lll | lll}
    & \multicolumn{3}{c|}{$p=0$} & \multicolumn{3}{c|}{$p=1$} & \multicolumn{3}{c|}{$p=2$}
    & \multicolumn{3}{c}{$p=3$}\\
    & $k_1$ & $k_2$ & $k_3$ & $k_1$ & $k_2$ & $k_3$ & $k_1$ & $k_2$ & $k_3$ & $k_1$ & $k_2$ & $k_3$\\
    \hline
    Hexagons &
        \highlightcell31 & \highlightcell53 & \highlightcell92 &
        \highlightcell25 & \highlightcell42 & \highlightcell80 &
        28 & \highlightcell47 & \highlightcell86 &
        28 & \highlightcell49 & \highlightcell90 \\
    Squares &
        37 & 64 & 116 &
        27 & 51 & 101 &
        \highlightcell27 & 51 & 98 & 
        \highlightcell27 & 52 & 100 \\
    Right triangles &
        40 & 70 & 134 &
        33 & 61 & 123 &
        31 & 60 & 117 &
        29 & 59 & 115\\
    Equilateral triangles &
        39 & 67 & 124 &
        33 & 58 & 113 & 
        32 & 59 & 113 &
        31 & 57 & 111 \\
    \end{tabular}
    \vspace{12pt}
    \caption{Iterations required for the GMRES iterative method with block 
    Jacobi preconditioner to converge. The smallest number of iterations in 
    each column is highlighted.}
    \label{tab:gmres-jac}
\end{table}

We now consider the solution of the above problem using the GMRES method with
the block ILU(0) preconditioner. Because of the sensitivity of the block ILU(0)
factorization to the ordering of the mesh elements, and for the sake of a fair
comparison between the generating patterns, we consider the \textit{natural
ordering} of mesh elements, illustrated in Figure \ref{fig:element-ordering}.
As in the case of the block Jacobi preconditioner, we repeat the test case of
the advection equation with variable velocity field. We record the number of
GMRES iterations required to converge to the above tolerance using the block ILU(0)
preconditioner in Table \ref{tab:gmres-ilu}. In this case, the square mesh
resulting in the smallest number of iterations in all of the trials. The mesh
consisting of right isosceles triangles resulted in the largest number of
iterations in all trials. We further note that the number of GMRES iterations 
required when using the block Jacobi preconditioner scales similarly to the 
number of block Jacobi iterations required, as recorded in Table 
\ref{tab:gaussian-iters}. We note that the block ILU(0) preconditioner 
requires fewer GMRES iterations to converge, and the number of iterations scales
more favorably in $k$, when compared with the block Jacobi 
preconditioner.

\begin{figure}[h]
    \centering
    \hspace*{\fill}%
    \begin{subfigure}{0.4\textwidth}
        \centering
        \begin{tikzpicture}
            \begin{axis}[
                hide axis,
                xmin=-0.35,
                xmax=0.75,
                axis equal,
                width=\textwidth,
                height=\textwidth,
                no markers,
                ticks=none,
            ]
            \addplot[black] table {hex.dat};
            \node at (axis cs:-0.14, 0) {1};
            \node at (axis cs:0.14,  0) {2};
            \node at (axis cs:0.42,  0) {3};
            
            \node at (axis cs:0,   0.25) {4};
            \node at (axis cs:0.28,0.25) {5};
            \node at (axis cs:0.56,0.25) {6};
            
            \node at (axis cs:-0.14, 0.5) {7};
            \node at (axis cs:0.14,  0.5) {8};
            \node at (axis cs:0.42,  0.5) {9};
            \end{axis}
        \end{tikzpicture}
        \caption{Hexagonal mesh}
    \end{subfigure}\hfill%
    \begin{subfigure}{0.4\textwidth}
        \centering
        \begin{tikzpicture}
            \begin{axis}[
                hide axis,
                xmin=-0.15,
                xmax=0.95,
                axis equal,
                width=\textwidth,
                height=\textwidth,
                no markers,
                ticks=none,
            ]
            \addplot[black] table {quad.dat};
            \node at (axis cs:0.13,  0.13) {1};
            \node at (axis cs:0.4,  0.13) {2};
            \node at (axis cs:0.65, 0.13) {3};
            
            \node at (axis cs:0.13,  0.4) {4};
            \node at (axis cs:0.4,  0.4) {5};
            \node at (axis cs:0.65, 0.4) {6};
            
            \node at (axis cs:0.13,  0.65) {7};
            \node at (axis cs:0.4,  0.65) {8};
            \node at (axis cs:0.65, 0.65) {9};
            \end{axis}
        \end{tikzpicture}
        \caption{Square mesh}
    \end{subfigure}%
    \hspace*{\fill}%
    
    \hspace*{\fill}%
    \begin{subfigure}{0.4\textwidth}
        \centering
        \begin{tikzpicture}
            \begin{axis}[
                hide axis,
                xmin=-0.15,
                xmax=0.95,
                axis equal,
                width=\textwidth,
                height=\textwidth,
                no markers,
                ticks=none,
            ]
            \addplot[black] table {rtri.dat};
            \node at (axis cs:0.125,  0.248) {1};
            \node at (axis cs:0.248,  0.125) {2};
            \node at (axis cs:0.496,  0.248) {3};
            \node at (axis cs:0.640,  0.125) {4};
            
            \node at (axis cs:0.125,  0.620) {5};
            \node at (axis cs:0.248,  0.496) {6};
            \node at (axis cs:0.496,  0.620) {7};
            \node at (axis cs:0.640,  0.496) {8};
            \end{axis}
        \end{tikzpicture}
        \caption{Right triangular mesh}
    \end{subfigure}\hfill%
    \begin{subfigure}{0.4\textwidth}
        \centering
        \begin{tikzpicture}
            \begin{axis}[
                hide axis,
                xmin=-0.3,
                xmax=0.8,
                axis equal,
                width=\textwidth,
                height=\textwidth,
                no markers,
                ticks=none,
            ]
            \addplot[black] table {etri.dat};
            \node at (axis cs:0,  0.231) {1};
            \node at (axis cs:0.2,  0.115) {2};
            \node at (axis cs:0.4,  0.231) {3};
            \node at (axis cs:0.6,  0.1155) {4};
            
            \node at (axis cs:0,  0.461) {5};
            \node at (axis cs:0.2,  0.577) {6};
            \node at (axis cs:0.4,  0.461) {7};
            \node at (axis cs:0.6,  0.577) {8};
            \end{axis}
        \end{tikzpicture}
        \caption{Equilateral triangular mesh}
    \end{subfigure}
    \hspace*{\fill}%
    
    \caption{Illustration of the natural ordering of mesh elements.}
    \label{fig:element-ordering}
\end{figure}
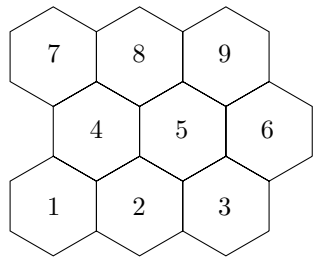
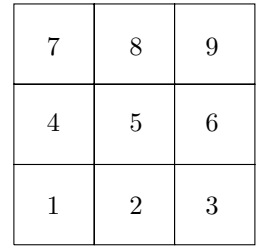
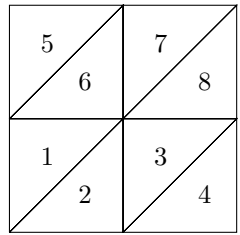
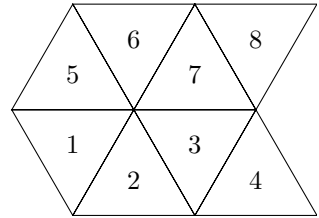

\begin{table}[h]
    \centering
    \tiny
    \setlength{\tabcolsep}{8pt}
    \begin{tabular}{r | lll | lll | lll | lll}
    & \multicolumn{3}{c|}{$p=0$} & \multicolumn{3}{c|}{$p=1$} & \multicolumn{3}{c|}{$p=2$}
    & \multicolumn{3}{c}{$p=3$}\\
    & $k_1$ & $k_2$ & $k_3$ & $k_1$ & $k_2$ & $k_3$ & $k_1$ & $k_2$ & $k_3$ & $k_1$ & $k_2$ & $k_3$\\
    \hline
    Hexagons &
        \highlightcell8 & 11 & \highlightcell16 &
        10 & 13 & 20 &
        11 & 15 & 23 &
        10 & 13 & 22 \\
    Squares &
        \highlightcell8 & \highlightcell10 & \highlightcell16 &
        \highlightcell8 & \highlightcell11 & \highlightcell19 &
        \highlightcell7 & \highlightcell10 & \highlightcell17 & 
        \highlightcell8 & \highlightcell10 & \highlightcell18 \\
    Right triangles &
        13 & 19 & 32 &
        10 & 14 & 28 &
        10 & 15 & 27 &
        11 & 14 & 28\\
    Equilateral triangles &
        11 & 15 & 27 &
        10 & 12 & 22 & 
        9 & 12 & 22 &
        9 & 12 & 22 \\
    \end{tabular}
    \vspace{12pt}
    \caption{Iterations required for the GMRES iterative method with ILU(0) preconditioner to converge. 
    The smallest number of iterations in each column is highlighted.}
    \label{tab:gmres-ilu}
\end{table}

\subsection{Compressible Euler equations}
The compressible Euler equations of gas dynamics in two dimensions (see \textit{e.g.} \cite{Hartmann}) 
are given by 
\begin{equation}
    \bm{u}_t + \nabla\cdot\bm{f}(\bm{u}) = 0,
\end{equation}
for
\begin{equation}
    \label{eq:euler-flux}
    \def\arraystretch{1}
    \bm{u} = \left(\begin{array}{c} \rho \\ \rho u \\ \rho v \\ \rho E\end{array}\right),
    \qquad
    \bm{f}_1(\bm{u}) = \left(\begin{array}{c} \rho u \\ \rho u^2 + p \\ \rho uv \\ \rho Hu\end{array}\right),
    \qquad
    \bm{f}_2(\bm{u}) = \left(\begin{array}{c} \rho v \\ \rho uv \\ \rho v^2 + p \\ \rho Hv\end{array}\right),
\end{equation}
where $\rho$ is the density, $\bm{v} = (u, v)$ is the fluid velocity, $p$ is the pressure, and $E$ is the specific 
energy. The total enthalpy $H$ is given by 
\begin{equation}
    H = E + \frac{p}{\rho},
\end{equation}
and the pressure is determined by the equation of state
\begin{equation}
    p = (\gamma - 1)\rho \left(E - \frac{1}{2}\bm{v}^2\right),
\end{equation}
where $\gamma = c_p/c_v$ is the ratio of specific heat capacities at constant pressure and constant volume.

We consider the model problem of an unsteady compressible vortex in a rectangular domain \cite{Wang2013}. 
The domain is taken to be a $20\times15$ rectangle and the vortex is initially centered at $(x_0, y_0) = (5, 5)$. 
The vortex is moving with the free-stream at an angle of $\theta$. The exact solution is given by
\begin{gather}
    u = u_\infty \left( \cos(\theta) - \frac{\epsilon ((y-y_0) - \overline{v} t)}{2\pi r_c} 
        \exp\left( \frac{f(x,y,t)}{2} \right) \right),\\
    u = u_\infty \left( \sin(\theta) - \frac{\epsilon ((x-x_0) - \overline{u} t)}{2\pi r_c} 
        \exp\left( \frac{f(x,y,t)}{2} \right) \right),\\
    \rho = \rho_\infty \left( 1 - \frac{\epsilon^2 (\gamma - 1)M^2_\infty}{8\pi^2} \exp((f(x,y,t))
        \right)^{\frac{1}{\gamma-1}}, \\
    p = p_\infty \left( 1 - \frac{\epsilon^2 (\gamma - 1)M^2_\infty}{8\pi^2} \exp((f(x,y,t))
        \right)^{\frac{\gamma}{\gamma-1}},
\end{gather}
where $f(x,y,t) = (1 - ((x-x_0) - \overline{u}t)^2 - ((y-y_0) - \overline{v}t)^2)/r_c^2$, 
$M_\infty$ is the Mach number, $u_\infty, \rho_\infty,$ and $p_\infty$ are the 
free-stream velocity, density, and pressure, respectively. The free-stream velocity is 
given by $(\overline{u}, \overline{v}) = u_\infty (\cos(\theta), \sin(\theta))$. The 
strength of the vortex is given by $\epsilon$, and its size is $r_c$. 
We choose the parameters to be $\gamma = 1.4$, $M_\infty = 0.5$, $u_\infty = 1$, $\theta = \arctan(1/2)$, 
$\epsilon = 0.3$, and $r_c = 1.5$. 

In the discontinuous Galerkin discretization of the Euler equations we use the Lax-Friedrichs numerical 
flux defined by
\begin{equation}
    \widehat{\bm{F}}(\bm{u}^+, \bm{u}^-, \bm{n}) = \tfrac{1}{2} \left(
        \bm{f}(\bm{u}^-)\cdot \bm{n} + \bm{f}(\bm{u}^+)\cdot \bm{n} + \alpha(\bm{u}^- - \bm{u}^+)\right),
\end{equation}
where $\alpha$ is the maximum absolute eigenvalue over $\bm{u}^-$ and $\bm{u}^+$ of the matrix $B(\bm{u}, \bm{n})$ 
defined by
\begin{equation}
    B(\bm{u}, \bm{n}) = J_{\bm{f}_1} n_1 + J_{\bm{f}_2} n_2,
\end{equation}
where $J_{\bm{f}_1}$ and $J_{\bm{f}_2}$ are the Jacobian matrices of the components of the numerical 
flux function $\bm{f}$ defined in equation \eqref{eq:euler-flux}.

We use the backward Euler time discretization, but remark that
\eqref{eq:semi-disc-dg} results in a nonlinear set of equations, which are
solved using Newton's method. Each iteration of Newton's method requires solving
a linear equation of the form \eqref{eq:be}. We set $h = 1$, and consider three
time steps, $k_1 = 0.03h$, $k_2 = 2 k_1$, $k_3 = 4 k_1$. We use piecewise
polynomials of degrees $p = 0, 1, 2, 3$. Each Newton solve requires between 3 to
8 iterations to converge to within a tolerance of 
$5 \times 10^{-13}$. The tolerance used for the linear solvers is the same as in 
the previous test cases.

\subsubsection{The Block Jacobi method}
Each iteration of Newton's method requires the solution of a linear system of 
equations. We solve these systems using the block Jacobi method. We compute
the  total the number of Jacobi iterations required to complete one solve of 
Newton's method, and report the results in Table \ref{tab:euler-jac}. We note 
that for each choice of $p$ and time step $k$, the hexagonal mesh required the 
fewest number of block Jacobi iterations. As in the previous numerical 
experiments, we do not see a decrease in performance for the hexagonal elements 
in the case of $p=3$. The square mesh resulted in the second-smallest number of 
iterations for most of the cases considered, while the two configurations of 
triangles resulted in generally similar numbers of iterations.

\begin{table}[h]
    \centering
    \tiny
    \setlength{\tabcolsep}{8pt}
    \begin{tabular}{r | lll | lll | lll | lll}
    & \multicolumn{3}{c|}{$p=0$} & \multicolumn{3}{c|}{$p=1$} & \multicolumn{3}{c|}{$p=2$}
    & \multicolumn{3}{c}{$p=3$}\\
    & $k_1$ & $k_2$ & $k_3$ & $k_1$ & $k_2$ & $k_3$ & $k_1$ & $k_2$ & $k_3$ & $k_1$ & $k_2$ & $k_3$\\
            \hline
            Hexagons                & 
        \highlightcell32 & \highlightcell49 & \highlightcell78 & 
        \highlightcell31 & \highlightcell50 & \highlightcell83 &
        \highlightcell50 & \highlightcell90 & \highlightcell158 & 
        \highlightcell53 & \highlightcell97 & \highlightcell171 \\
            Squares                 & 
        34 & 51 & 89 &
        \highlightcell31 & 54 & 92 &
        54 & 99 & 181 &
        55 & 105 & 201  \\
            Right triangles         & 
        37 & 56 & 97 &
        41 & 64 & 112 &
        58 & 101 & 189 &
        59 & 113 & 217 \\
            Equilateral triangles   & 
        37 & 57 & 95 &
        39 & 62 & 113 &
        54 & 99 & 179 &
        60 & 114 & 215
        \end{tabular}
        \vspace{8pt}
    \caption{Block Jacobi iterations required per Newton solve of the 
    compressible Euler equations. The lowest number of iterations in each 
    column is highlighted.}
    \label{tab:euler-jac}
\end{table}

\subsubsection{The GMRES method}
We now repeat the above test case, using the GMRES method to solve the resulting
linear systems. We consider both the block Jacobi and block ILU(0)
preconditioners. We then compute the total number of GMRES iterations required
to complete one solve of Newton's method. As in Section \ref{sec:adv-gmres}, the
ordering of the mesh elements has a significant effect on the effectiveness of
the block ILU(0) approximate factorization. For this reason, we use the natural
ordering of elements, depicted in Figure \ref{fig:element-ordering}. We present
the results for the block Jacobi preconditioner in Table
\ref{tab:euler-gmres-jac}, and for the block ILU(0) preconditioner in Table
\ref{tab:euler-gmres-ilu}. With the block Jacobi preconditioner, the hexagonal
mesh required the smallest number of iterations for all test cases considered,
and the square mesh the second-smallest. In the case of the block ILU(0)
preconditioner, the square mesh required the fewest number of iterations, with
the hexagonal mesh usually requiring the second-smallest number of iterations.
As we observed in Section \ref{sec:adv-gmres}, the number of iterations required
for both the block Jacobi method and GMRES with the block Jacobi preconditioner
scales quite poorly with increasing timesteps. The number of GMRES iterations
required when using the block ILU(0) preconditioner is significantly better.

\begin{table}[h]
    \centering
    \tiny
    \setlength{\tabcolsep}{8pt}
    \begin{tabular}{r | lll | lll | lll | lll}
    & \multicolumn{3}{c|}{$p=0$} & \multicolumn{3}{c|}{$p=1$} & \multicolumn{3}{c|}{$p=2$}
    & \multicolumn{3}{c}{$p=3$}\\
    & $k_1$ & $k_2$ & $k_3$ & $k_1$ & $k_2$ & $k_3$ & $k_1$ & $k_2$ & $k_3$ & $k_1$ & $k_2$ & $k_3$\\
            \hline
            Hexagons                & 
        \highlightcell55 & \highlightcell74 & \highlightcell106 & 
        \highlightcell50 & \highlightcell92 & \highlightcell126 &
        \highlightcell61 & \highlightcell110 & \highlightcell153 & 
        \highlightcell76 & \highlightcell141 & \highlightcell195 \\
            Squares                 & 
        62 & 84 & 155 &
        52 & 93 & 132 &
        67 & 126 & 185 &
        78 & 149 & 222 \\
            Right triangles         & 
        63 & 87 & 162 &
        81 & 106 & 184 &
        96 & 132 & 242 &
        85 & 159 & 299\\
            Equilateral triangles   & 
        66 & 90 & 167 &
        81 & 108 & 187 &
        72 & 133 & 197 &
        85 & 161 & 245
        \end{tabular}
        \vspace{8pt}
    \caption{GMRES with block Jacobi preconditioner. Iterations required per 
    Newton solve of the compressible Euler equations. The lowest number of 
    iterations in each column is highlighted.}
    \label{tab:euler-gmres-jac}
\end{table}

\begin{table}[h]
    \centering
    \tiny
    \setlength{\tabcolsep}{8pt}
    \begin{tabular}{r | lll | lll | lll | lll}
    & \multicolumn{3}{c|}{$p=0$} & \multicolumn{3}{c|}{$p=1$} & \multicolumn{3}{c|}{$p=2$}
    & \multicolumn{3}{c}{$p=3$}\\
    & $k_1$ & $k_2$ & $k_3$ & $k_1$ & $k_2$ & $k_3$ & $k_1$ & $k_2$ & $k_3$ & $k_1$ & $k_2$ & $k_3$\\
    \hline
            Hexagons                & 
        \highlightcell24 & 32 & \highlightcell42 & 
        \highlightcell21 & 36 & 48  &
        29 & 48 & 57 &
        29 & 50 & 64 \\
            Squares                 & 
        \highlightcell24 & \highlightcell28 & 45 &
        \highlightcell21 & \highlightcell33 & \highlightcell40 &
        \highlightcell24 & \highlightcell41 & \highlightcell49 &
        \highlightcell27 & \highlightcell48 & \highlightcell60 \\
            Right triangles         & 
        31 & 40 & 70 &
        35 & 40 & 60 &
        36 & 48 & 69 &
        31 & 49 & 75 \\
            Equilateral triangles   & 
        28 & 37 & 65 &
        37 & 44 & 70 &
        33 & 56 & 68 &
        38 & 64 & 80  \\
        \end{tabular}
        \vspace{8pt}
    \caption{GMRES with block ILU(0) preconditioner. Iterations required per 
    Newton solve of the compressible Euler equations. The lowest number of 
    iterations in each column is highlighted.}
    \label{tab:euler-gmres-ilu}
\end{table}

\subsection{Inviscid flow problems}

The following two numerical experiments extend the above results to
larger-scale, more realistic flow problems. These problems, in contrast to the
preceding test cases, are characterized by a large number of degrees of freedom,
the presence of geometric features and wall boundary conditions, variably-sized
mesh elements, and shocks. As in the previous section, the equations considered
here are the compressible Euler equations. For the following two problems, we
choose the finite element function space to consist of piecewise constant
functions (corresponding to $p=0$), which results in a finite-volume-type
discretization. This choice of discretization allows for the solution of
problems with shocks, without the use of slope limiters, artificial viscosity,
or other shock-capturing techniques \cite{leveque2002finite}. The Roe numerical
flux is used as an approximate Riemann solver for these problems.

\subsubsection{Subsonic flow over a circular cylinder}

For a first test case, we consider the inviscid flow over a circular cylinder at
Mach 0.2. The computational domain is defined as $\Omega = R \setminus C$, where
$R = [-10, 30] \times [-10, 20]$, and $C$ is a disk of radius 1 centered at the
point $(5,5)$. Farfield boundary conditions are enforced on $\partial R$, and a
no normal flow condition is enforced on $\partial C$. The freestream velocity is
taken to be unity in the $x$-direction, and $\rho_\infty = 1$. For this test
case we use four unstructured meshes, two consisting entirely of triangles, and
two consisting of mixed polygons, generated using the PolyMesher algorithm
\cite{Talischi:2012}. All the meshes are created using a gradient-limited
element size function that determines the initial distribution of seed points
according to the rejection method \cite{persson2005mesh}, such that the element
edge length near the surface of the cylinder is about one-fifth the edge length
of elements away from the cylinder. For both the triangular and polygonal
meshes, we consider a coarse mesh, with 15{,}404 elements, and a fine mesh with
62{,}270 elements. Thus, the average area of each element is the same for both
the polygonal and triangular meshes. Additionally, the number of degrees of
freedom in the solution is the same, allowing for a fair comparison. The coarse
polygonal mesh, and a zoom-in around the surface of the cylinder are shown in
Figure \ref{fig:cyl-mesh}.

\begin{figure}[t]
    \centering
    \begin{subfigure}{3.35in}
        \centering
        \includegraphics[width=3.3in]{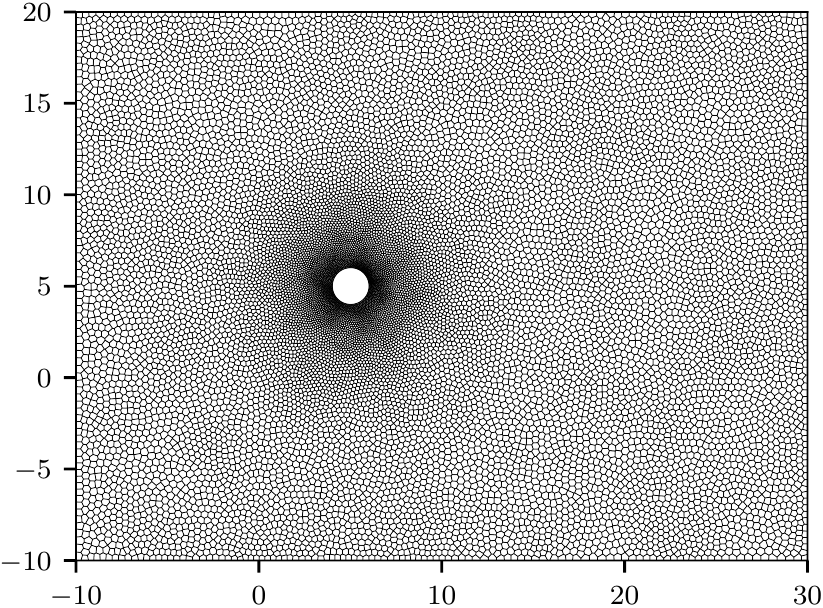}
    \end{subfigure}\hfill%
    \begin{subfigure}{2.55in}
        \centering
        \includegraphics[width=2.5in]{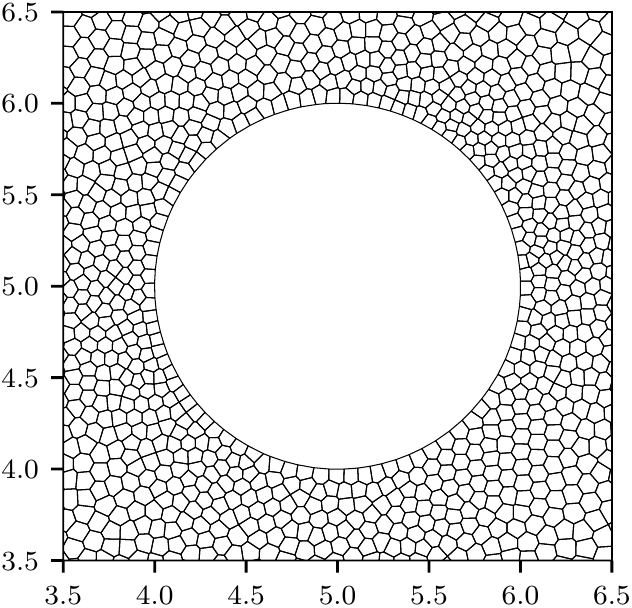}
    \end{subfigure}%
    
    \caption{Overview of the coarse mesh with 15{,}404 elements, with zoom-in 
             showing polygonal elements near the surface of the cylinder.}
    \label{fig:cyl-mesh}
\end{figure}

Starting from freestream initial conditions, we integrate the equations until $t
= 5 \times 10^{-3}$ in order to obtain a representative solution. Using this
solution, we then compute 10 time steps using a third-order $A$-stable DIRK
method \cite{Alexander1977}. Each stage of the DIRK method requires the solution
of a nonlinear system of equations, which we solve by means of Newton's method.
In each iteration of Newton's method, we solve the resulting linear system of
the form \eqref{eq:be} using both the block Jacobi method and the preconditioned
GMRES method.
The nonlinear system is solved to within a tolerance of 
$10^{-8}$, and each linear system is solved using a relative tolerance of $10^{-5}$.
 For the GMRES method, we consider two preconditioners: block
Jacobi, and block ILU(0). In order to compare the iterative solver performance
differences between meshes, we compute the total number of solver iterations
required to complete all 10 time steps. The results for the GMRES method are 
shown in Table \ref{tab:cylinder-gmres}, and for the block Jacobi solver in
Table \ref{tab:cylinder-jacobi}.

\begin{table}[t] 
\caption{Total GMRES iterations per 10 time steps for inviscid flow
         over a circular cylinder.}
\label{tab:cylinder-gmres}
\centering
\begin{subtable}{0.99\linewidth}
\centering
\caption{Coarse grid with 15{,}404 elements}
\begin{tabular}{lcccc|cc}
\toprule
 & \multicolumn{2}{c}{ILU} & \multicolumn{2}{c|}{Jacobi} & 
  \multicolumn{2}{c}{Ratios} \\
$\Delta t$ & Polygonal & Triangular & Polygonal & Triangular & ILU & Jacobi \\
\midrule
$1.0 \times 10^{-1}$ & 793 & 932 & 2092 & 3126 & 0.85 & 0.67 \\
$2.5 \times 10^{-1}$ & 1569 & 1829 & 4405 & 6870 & 0.86 & 0.64 \\
$5.0 \times 10^{-1}$ & 2470 & 3090 & 7145 & 11859 & 0.80 & 0.60 \\
$1.0$                & 3651 & 4486 & 11054 & 18880 & 0.81 & 0.59 \\
\bottomrule
\end{tabular}
\end{subtable} \\
\begin{subtable}{0.99\linewidth}
\centering
\caption{Fine mesh with 95{,}932 elements}
\begin{tabular}{lcccc|cc}
\toprule
 & \multicolumn{2}{c}{ILU} & \multicolumn{2}{c|}{Jacobi} & 
  \multicolumn{2}{c}{Ratios} \\
$\Delta t$ & Polygonal & Triangular & Polygonal & Triangular & ILU & Jacobi \\
\midrule
$1.0 \times 10^{-1}$ & 1443 & 1673 & 4075 & 6137 & 0.86 & 0.66 \\
$2.5 \times 10^{-1}$ & 2998 & 3344 & 8732 & 12741 & 0.90 & 0.69 \\
$5.0 \times 10^{-1}$ & 4720 & 5423 & 14084 & 21882 & 0.87 & 0.64 \\
$1.0$                & 7205 & 8151 & 22814 & 34706 & 0.88 & 0.66 \\
\bottomrule
\end{tabular}
\end{subtable}
\end{table}
\begin{table}[t] 
\caption{Total block Jacobi iterations per 10 time steps for inviscid 
         flow over a circular cylinder.}
\label{tab:cylinder-jacobi}
\centering
\begin{subtable}{0.49\linewidth}
\centering
\caption{Coarse grid with 15{,}404 elements}
\begin{tabular}{llll}
\toprule
$\Delta t$ & Polygonal & Triangular & Ratio \\
\midrule
$1.0 \times 10^{-1}$ & 2474 & 3159 & 0.78 \\
$2.5 \times 10^{-1}$ & 4895 & 6697 & 0.73 \\
$5.0 \times 10^{-1}$ & 7882 & 12158 & 0.65 \\
$1.0$              & 13181 & 19072 & 0.69 \\
\bottomrule
\end{tabular}
\end{subtable}
\begin{subtable}{0.49\linewidth}
\centering
\caption{Fine mesh with 95{,}932 elements}
\begin{tabular}{llll}
\toprule
$\Delta t$ & Polygonal & Triangular & Ratio \\
\midrule
$1.0 \times 10^{-1}$ & 4788 & 6281 & 0.76 \\
$2.5 \times 10^{-1}$ & 9609 & 12406 & 0.77 \\
$5.0 \times 10^{-1}$ & 15580 & 20946 & 0.74 \\
$1.0$                & 26628 & 33934 & 0.78 \\
\bottomrule
\end{tabular}
\end{subtable}
\end{table}

These results demonstrate a consistent trend, corroborating both the numerical
results and the analysis from the previous sections. When using the block Jacobi
solver or GMRES with block Jacobi preconditioner, the polygonal mesh results in
convergence in between 60--70\% of the iterations required for the triangular
mesh. The effect is smaller when using the ILU(0) preconditioner, but we do
still observe a modest reduction in the number of iterations required. When 
using the block Jacobi iterative solver, we observe iteration counts very 
similar to when using GMRES with block Jacobi as a preconditioner. In these 
cases, the polygonal mesh requires between 70--80\% of the iterations as the 
all-triangular mesh.

\subsubsection{Supersonic flow over a circular cylinder}

The next numerical example is designed to investigate the performance of the
iterative solvers for steady-state problems, in the presence of shocks and
$h$-adapted meshes. For this problem, we let the domain be $\Omega = R \setminus
C$, where $ R =[0,5]\times[0,10]$, and, as before, $C$ is a circle of radius one
centered at $(5,5)$. Freestream conditions are enforced at the left, top, and
bottom boundaries, an inviscid wall condition is enforced on the boundary of the
cylinder, and an outflow condition is enforced on the right boundary. The Mach
number is set to $M = 2.0$, resulting in the formation of a shock upstream from
the cylinder. In order to accurately capture the shock, we refine the mesh in
its vicinity. As in the previous case, we consider a set of four meshes, two
all-triangular, and two polygonal. For both the triangular and polygonal meshes,
we consider coarse and fine versions, with 31{,}162 and 95{,}932 elements,
respectively. The coarse mesh is depicted in Figure \ref{fig:supersonic-mesh},
with Mach isolines overlaid to indicate the position of the shock. Additionally,
Mach contours of the steady-state solution are shown in Figure
\ref{fig:supersonic-contours}.

\begin{figure}[t]
    \centering
    \begin{subfigure}{0.49\textwidth}
        \centering
        \includegraphics[width=1.7in]{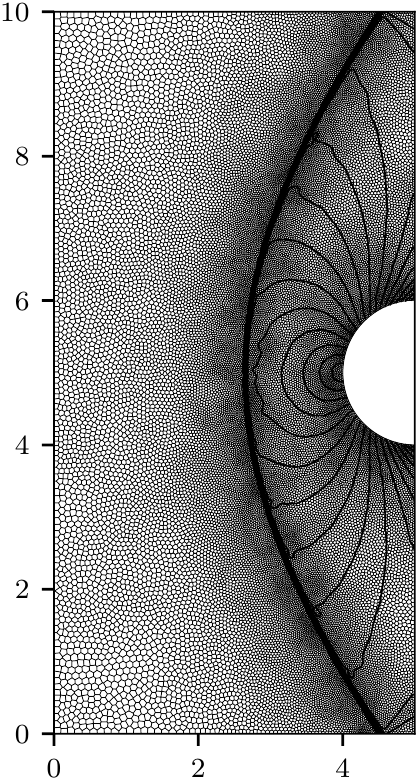}
        \caption{Coarse mesh for supersonic test problem, showing Mach isolines
                 for steady-state solution}
        \label{fig:supersonic-mesh}
    \end{subfigure}\hfill%
    \begin{subfigure}{0.49\textwidth}
        \centering
        \includegraphics[width=2.3in]{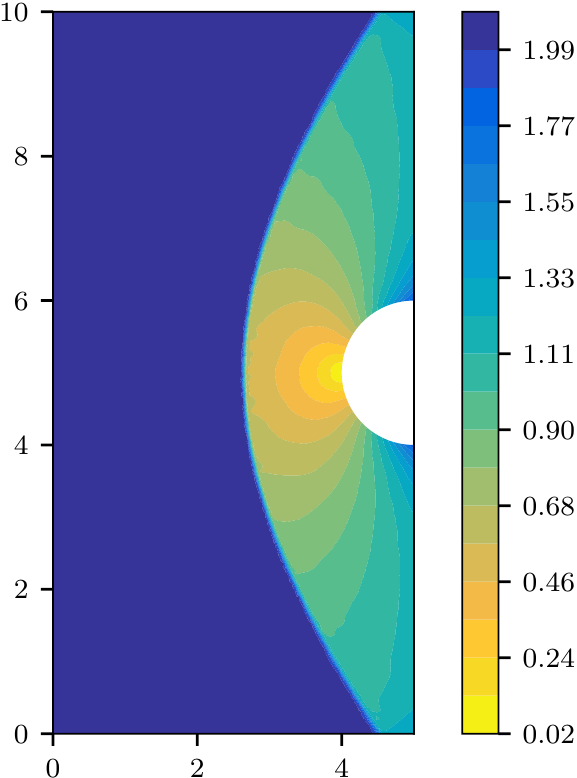}
        \caption{Contours of Mach number for steady state solution}
        \label{fig:supersonic-contours}
    \end{subfigure}%
    \caption{Overview of coarse polygonal mesh with 31{,}162 elements, showing 
             Mach number contours for steady-state solution.}
    \label{fig:supersonic}
\end{figure}

Beginning with freestream initial conditions, the solution rapidly approaches a
steady state. We integrate in time until $t = 100$ in order to obtain an
solution which can be used as an initial guess for the steady-state Newton
solve. Then, starting with this solution, we set the time-derivative of the
solution to zero and solve the resulting nonlinear equations using Newton's
method to find a steady-state solution. The resulting linear system that is
required to be solved at each iteration can be thought of as corresponding to
equation \eqref{eq:be}, where formally we set $k = \infty$.
The nonlinear system is solved to within a tolerance of 
$10^{-10}$, and each linear system is solved using a relative tolerance of $10^{-5}$.
Since the mass
matrix in \eqref{eq:be} acts to regularize the linear system, the conditioning
becomes worse for larger values of $k$, and the number of iterations required
per linear solve grows. Hence, effective preconditioners are particularly
important for the solution of such steady-state problems. For these problems, 
the block Jacobi iterative solver did not converge in fewer than 10{,}000 
iterations, and so we consider only the GMRES method, using block ILU(0) and 
block Jacobi preconditioners.

We present the comparison of iteration counts for this problem in Table 
\ref{tab:supersonic-results}. On the coarse meshes, the ILU(0) preconditioner 
required about 73\% as many iterations on the polygonal mesh when compared with
the triangular mesh. This difference is more significant when using the block 
Jacobi preconditioner, consistent with the results observed in previous section.
In this case, the polygonal mesh requires only slightly more than one third the 
number of iterations as the all-triangular mesh. On the fine mesh, there are 
close to half a million degrees of freedom. For a problem of this scale, we did 
not observe convergence in less than 10{,}000 iterations per linear solve using 
the block Jacobi preconditioner, and so we only compare performance using the 
block ILU(0) preconditioner. In this case, the polygonal mesh required about 
half as many iterations per steady-state solve when compared with the 
all-triangular mesh.

\begin{table}[h!] 
\caption{Total GMRES iterations per steady-state solve for supersonic 
         flow over a cylinder.}
\label{tab:supersonic-results}
\centering
\begin{subtable}{0.49\linewidth}
\caption{Coarse grid with 31{,}162 elements}
\begin{tabular}{lccc}
\toprule
 & Polygonal & Triangular & Ratio \\
\midrule
ILU & 469 & 640 & 0.73 \\
Jacobi & 2340 & 6464 & 0.36 \\
\bottomrule
\end{tabular}
\end{subtable}
\begin{subtable}{0.49\linewidth}
\caption{Fine mesh with 95{,}932 elements}
\begin{tabular}{lccc}
\toprule
 & Polygonal & Triangular & Ratio \\
\midrule
ILU & 953 & 1947 & 0.49 \\
Jacobi & -- & -- & -- \\
\bottomrule
\end{tabular}
\end{subtable}
\end{table}
}

\section{Conclusions}\label{sec:conclusions}
In this paper we have analyzed the effect of the generating pattern of a regular
mesh on the convergence of iterative linear solvers applied to implicit
discontinuous Galerkin discretizations. We considered four generating patters: a
hexagon, a square, two right triangles, and two equilateral triangles.

A classical von Neumann analysis applied to the constant-velocity advection
equation allowed us to compute the eigenvalues of the block Jacobi matrix, and
therefore estimate the speed of convergence of the block Jacobi method. In more
than half of the cases considered, the hexagonal generating pattern resulted in
the smallest eigenvalues, and in the remaining cases, the square generating
pattern resulted in the smallest eigenvalues.

In order to extend these results beyond the case of the constant-velocity
advection equation, we performed numerical experiments on the variable-velocity
advection equation and compressible Euler equations. In the case of the
advection equation, in all but one case the hexagonal mesh resulted in the
fastest convergence, and in the remaining case the square mesh resulted in the
fastest convergence. In the case of the Euler equations, the hexagonal mesh
resulted in the fastest convergence in all test cases.

We additionally considered two irregular meshes resulting from the random
perturbation of a set of regularly-spaced generating points. We obtain a
triangular mesh by performing the Delaunay triangulation on these points, and we
obtain a polygonal mesh by constructing the Voronoi diagram dual to the Delaunay
triangulation. Solving the advection equation on these irregular meshes, we
observed that the block Jacobi method converged faster on the polygonal mesh in
every test case. Additionally, we performed numerical experiments examining the
performance of the GMRES iterative method when used with the ILU(0)
preconditioner. We found that in all of the test cases, the square generating
pattern resulted in the fewest number of GMRES iterations, and in all but two
cases, the hexagonal generating pattern resulted in the second-fewest number of
iterations.

For a final set of numerical experiments, we performed two inviscid
fluid flow simulations on sets of coarse and fine meshes. Each mesh was either
all-triangular, or was composed of arbitrary polygons. We measured iteration
counts for both time-dependent and steady-state problems, using the block Jacobi
method, and GMRES with block ILU(0) and block Jacobi preconditioners. We found
that the polygonal meshes resulted in faster convergence of the iterative
solvers, with a larger difference being observed for the block Jacobi method and
preconditioner. This difference was more pronounced for the steady-state
problem, with a quite significant difference observed on the fine mesh using
GMRES with ILU(0).

These results suggest that certain types of polygonal meshes have the advantage
of rapid convergence of iterative solvers. Future research directions involve
the study of accuracy of DG methods on polygonal and polyhedral meshes,
efficient computation of quadrature rules over arbitrary polygonal domains and
the extension of the above results to three spatial dimensions.

\bibliographystyle{plain}
\bibliography{polygonal}

\end{document}